\documentclass[12pt]{article}
\usepackage{booktabs}
\usepackage{caption}
\usepackage{mathrsfs}
\usepackage{amsmath}
\usepackage{amsfonts,amsthm,amssymb,mathrsfs,bbding}
\usepackage{graphics,multicol}
\usepackage{graphicx}
\usepackage{color}
\usepackage{enumerate}
\usepackage{caption}
\usepackage{rotating}
\usepackage{lscape}
\usepackage{longtable}
\allowdisplaybreaks[4]
\usepackage{tabularx}
\usepackage[%dvipdfm,
colorlinks,
anchorcolor=blue,
linkcolor=blue,
citecolor=blue]{hyperref}
\usepackage{extarrows}
\usepackage{cite}
\usepackage{latexsym,bm}
\usepackage{mathtools}
\evensidemargin=\oddsidemargin\topmargin=-1.5cm

\newtheorem{thm}{Theorem}

\newtheorem{lem}[thm]{Lemma}

\newtheorem{remark}{Remark}

\newtheorem{claim}{Claim}
\theoremstyle{definition}
\renewcommand\proofname{\bf Proof}
\addtocounter{section}{0}
\begin{document}

\title{\bf The $A_\alpha$-spectral radius and perfect matchings of graphs}
\author{{Yanhua Zhao, Xueyi Huang\footnote{Corresponding author.}\setcounter{footnote}{-1}\footnote{\emph{E-mail address:} huangxymath@163.com.} \ and Zhiwen Wang}\\[2mm]
\small Department of Mathematics, East China University of Science and Technology, \\
\small  Shanghai 200237, P.R. China}

\date{}
\maketitle
{\flushleft\large\bf Abstract} Let $\alpha\in[0,1)$, and let $G$ be a graph of even order $n$ with  $n\geq f(\alpha)$, where $f(\alpha)=10$ for  $0\leq \alpha\leq1/2$, $f(\alpha)=14$ for $1/2<\alpha\leq 2/3$ and $f(\alpha)=5/(1-\alpha)$ for $2/3<\alpha<1$. In this paper, it is shown that  if the $A_\alpha$-spectral radius of $G$ is not less than the largest root of $x^3 - ((\alpha + 1)n +\alpha-4)x^2 + (\alpha n^2 + (\alpha^2 - 2\alpha - 1)n - 2\alpha+1)x -\alpha^2n^2 + (5\alpha^2 - 3\alpha + 2)n - 10\alpha^2 + 15\alpha - 8=0$ then $G$ has a perfect matching unless $G=K_1\nabla(K_{n-3}\cup 2K_1)$. This generalizes a result of  S. O [Spectral radius and matchings in graphs, Linear Algebra Appl. 614 (2021) 316--324], which gives a sufficient condition for the existence of a perfect matching in a graph in terms of the adjacency spectral radius.
\begin{flushleft}
\textbf{Keywords:}  Perfect matching, $A_\alpha$-spectral radius,  quotient matrix.
\end{flushleft}
\textbf{AMS Classification:} 05C50

\section{Introduction}
Let $G$ be an undirected simple graph with vertex set $V(G)$ and edge set $E(G)$. The \textit{adjacency matrix} of $G$ is defined as $A(G)=(a_{u,v})_{u,v\in V(G)}$, where $a_{u,v}=1$ if $uv\in E(G)$, and $a_{u,v}=0$ otherwise. Let $D(G)$ denote the diagonal matrix of vertex degrees  of $G$. Then $L(G)=D(G)-A(G)$ and $Q(G)=D(G)+A(G)$ are called the \textit{Laplacian matrix}  and \textit{signless Laplacian matrix} of $G$, respectively. For any $\alpha \in [0,1)$, Nikiforov \cite{ni} introduced the $A_\alpha$-\textit{matrix} of $G$ as
$$A_\alpha(G)=\alpha D(G)+(1-\alpha)A(G).$$
Notice that $A_0(G)=A(G)$ and $A_{\frac{1}{2}}(G)=\frac{1}{2}Q(G)$. The eigenvalues of $A_\alpha(G)$ are called the $A_\alpha$-\textit{eigenvalues} of $G$, and the largest one of them, denoted by $\rho_\alpha(G)$,  is called the $A_\alpha$-\textit{spectral radius} of $G$.  The adjacency, Laplacian and signless Laplacian eigenvalues (and spectral radius) can be similarly defined.  For some interesting spectral properties of  $A_\alpha(G)$, we refer the reader to \cite{Lin1,Lin2,Lin3,Liu,ni, ni2, ni3,ni4}.

%For any $S\subseteq V(G)$, we denote by $G[S]$ the subgraph of $G$ induced by $S$, and $G-S$ the graph obtained from $G$ by deleting $S$ together with those edges incident to $S$. The \textit{complement graph}  $\overline{G}$ of $G$ is the graph with the same vertex set as $G$ in which two vertices are adjacent if and only if they are not adjacent in $G$. Let $G\cup H$ denote the disjoint union of two graphs $G$ and $H$. The \textit{join} of $G$ and $H$, denoted by $G\nabla H$, is the graph obtained from $G\cup H$ by adding all edges between $G$ and $H$. Also, we denote by $K_n$ the complete
%graph on $n$ vertices.

For any $S\subseteq V(G)$, we denote by  $G-S$ the graph obtained from $G$ by deleting $S$ together with those edges incident to $S$. Let $G\cup H$ denote the disjoint union of two graphs $G$ and $H$. The \textit{join} of $G$ and $H$, denoted by $G\nabla H$, is the graph obtained from $G\cup H$ by adding all edges between $G$ and $H$. Also, we denote by $K_n$ the complete graph on $n$ vertices. A \textit{matching} $M$ of  $G$ is a subset of $E(G)$ such that any two edges of $M$ have no common vertices. If each vertex of $G$ is incident to exactly one edge of $M$, then $M$ is called a \textit{perfect matching} of $G$. The \textit{matching number} $\alpha'(G)$ of $G$ is the number of edges in a maximum matching of $G$.

The relationships between the matching number and the eigenvalues of various graph matrices have been investigated by several researchers in the past two decades. In 2005, Brouwer and Haemers \cite{br} provided some sufficient conditions for the existence of a perfect matching in a graph in terms of the largest and the second smallest Laplacian eigenvalues.  Also, they proved that a $k$-regular graph $G$ of even order has a perfect matching if the third largest adjacency eigenvalue of $G$ is at most $k-1+3/(k+1)$ for $k$ even and $k-1+3/(k+2)$ for $k$ odd. Later, Cioab\u{a}, Gregory and Haemers \cite{ci,ci1,ci2}   improved and extended  this result, and gave some sufficient conditions for the existence of large matchings in regular graphs. In addition, the relationships between the fractional matching number (see \cite{SU} for the definition) and the adjacency, Laplacian and signless Laplacian spectral radius were studied in \cite{so1}, \cite{XZS} and \cite{py}, respectively.

%A \textit{fractional matching} of a graph $G$ is a function $f:E(G)\rightarrow [0,1]$ such that $\sum_{e\in N(v)}f(e) \leq 1$ for all $v \in V(G)$, where $N(v)$ is the set of neighbors of $v$. The \textit{fractional matching number} of $G$, denoted by $\alpha'_*(G)$, is the maximum value
%of $\sum_{e\in E(G)}f(e)$ over all fractional matchings. A \textit{fractional perfect matching} of a graph $G$ is a fractional matching $f$ with $\alpha'_*(G)= \sum_{e\in E(G)}f(e)=n/2$.

Very recently,  O \cite{so}  provided a lower bound for the adjacency spectral radius $\rho(G)$ which guarantees the existence of a perfect matching in a connected graph $G$.

\begin{thm}(O \cite{so})\label{thm-1}
Let $n \geq 8$ be an even integer or $n=4$. If $G$ is an $n$-vertex connected graph with
$\rho(G) >\theta(n)$, where $\theta(n)$ is the largest root of $x^3 -(n-4)x^2-(n-1)x+2(n-4)=0$, then
$G$ has a perfect matching. For $n = 6$, if $\rho(G) > \frac{1+\sqrt{33}}{2}$, then $G$ has a perfect matching.
\end{thm}

Liu, Pan and Li \cite{LPL} gave an analogue of Theorem \ref{thm-1} for the the signless Laplacian spectral radius $\rho_Q(G)$ of a connected graph $G$.
\begin{thm}(Liu, Pan and Li \cite{LPL})\label{thm-1-1}
Let $G$ be an $n$-vertex connected graph, where $n$ is an even number. Then
\begin{enumerate}[(i)]
\item $G$ has a perfect matching for $n\geq 10$ or $n=4$, if $\rho_Q(G)>r(n)$, where $r(n)$ is the largest root of equation $x^3 - (3n - 7)x^2 + n(2n - 7)x - 2(n^2
 - 7n + 12) = 0$;
\item $G$ has a perfect matching for $n = 6$ or $n = 8$, if  $\rho_Q(G)>4 + 2\sqrt{3}$ or $\rho_Q(G)>6+2\sqrt{6}$,
respectively.
\end{enumerate}
\end{thm}

It is natural to ask  whether the above results can be extended to the $A_\alpha$-spectral radius of a graph. In this paper, by using a similar method as in \cite{so}, we  prove the following result.

\begin{thm}\label{thm-2}
Let $\alpha\in[0,1)$, and let $G$ be a connected graph of even order $n$ with  $n\geq f(\alpha)$, where 
\begin{equation}\label{equ-0}
f(\alpha)=\left\{
\begin{array}{ll}
10,&\mbox{if $0\leq \alpha\leq1/2$};\\
14,&\mbox{if $1/2<\alpha\leq 2/3$};\\
5/(1-\alpha),&\mbox{if $2/3<\alpha<1$}.
\end{array}
\right.
\end{equation}
If $\rho_\alpha(G) \geq \rho_\alpha(K_1\nabla(K_{n-3}\cup 2K_1))$, then $G$ has a perfect matching unless $G=K_1\nabla(K_{n-3}\cup 2K_1)$, where $\rho_\alpha(K_1\nabla(K_{n-3}\cup 2K_1))$ is equal to the largest root of $x^3 - ((\alpha + 1)n +\alpha-4)x^2 + (\alpha n^2 + (\alpha^2 - 2\alpha - 1)n - 2\alpha+1)x -\alpha^2n^2 + (5\alpha^2 - 3\alpha + 2)n - 10\alpha^2 + 15\alpha - 8=0$.
\end{thm}

%It is worth mentioning that the extremal graph in Theorem \ref{thm-2} is exactly same as the one in \cite{so}. Furthermore,  from Theorem \ref{thm-2} we immediately deduce the following result.
%\begin{cor}\label{cor1}
%Let $G$ be a graph of even order $n\geq 16$, and let $\tau(n)$ be the largest root of $x^3 - (3n-7)/2\cdot x^2 + (2n^2 - 7n)/4\cdot x - (n^2-7n+12)/4=0$.  If the signless Laplacian radius of $G$ is not less than $2\cdot \tau(n)$, then
%$G$ has a perfect matching unless $G=K_1\nabla(K_{n-3}\cup 2K_1)$.
%\end{cor}

\begin{remark}
\emph{
 Recall that $A_0(G)=A(G)$ and $A_{\frac{1}{2}}(G)=\frac{1}{2}Q(G)$. For $\alpha=0$ and $\alpha=1/2$, one can verify that the conclusion of Theorem \ref{thm-2} coincides with the results of Theorem \ref{thm-1}  and Theorem \ref{thm-1-1}, respectively, provided that $n\geq 10$. 
}
\end{remark}

\section{Preliminaries}
In this section, we list three lemmas, which are useful in the proof of Theorem \ref{thm-2}.

Let $M$ be a real $n\times n$ matrix, and let $X=\{1,2,\ldots,n\}$. Given a partition $\Pi:X=X_1\cup X_2\cup \cdots \cup X_k$,  the matrix $M$ can be correspondingly partitioned as
$$
M=\left(\begin{array}{ccccccc}
M_{1,1}&M_{1,2}&\cdots &M_{1,k}\\
M_{2,1}&M_{2,2}&\cdots &M_{2,k}\\
\vdots& \vdots& \ddots& \vdots\\
M_{k,1}&M_{k,2}&\cdots &M_{k,k}\\
\end{array}\right).
$$
The \textit{quotient matrix} of $M$ with respect to $\Pi$, written $B_\Pi=(b_{i,j})_{i,j=1}^k$,  is defined by
$$b_{i,j}=\frac{1}{|X_i|}\mathbf{j}_{|X_i|}^TM_{i,j}\mathbf{j}_{|X_j|}$$
for all $i,j\in \{1,2,\ldots,k\}$, where $\mathbf{j}_{s}$ denotes the all ones vector in $\mathbb{R}^s$.
The partition $\Pi$ is called \textit{equitable} if each block $M_{i,j}$ of $M$ has constant row sum $b_{i,j}$.
Also, the quotient matrix $B_\Pi$ is called \textit{equitable} if $\Pi$ is an  equitable partition of $M$.

\begin{lem}(See \cite[Theorem 2.5]{yl})\label{lem-1}
Let $M$ be a nonnegative matrix, and let $B$ be an equitable quotient matrix of $M$. Then the eigenvalues of $B$ are also eigenvalues of $M$, and
$$\rho(M) = \rho(B),$$
where $\rho(M)$ and $\rho(B)$ are the spectral radius of $M$ and $B$, respectively.
\end{lem}
%Let $x_1\geq x_2 \geq ... \geq x_n$ and $y_1 \geq y_2 \geq ... \geq y_m$ be two sequences of real numbers with $n> m$. The second sequence
%is said to interlace the first one whenever $x_i\geq y_i \geq x_{n-m+i}$ for $i\in[1,m]$.
%\begin{lem}(\cite{hh})\label{lem-2}
% Let $B$ be a quotient matrix of a symmetric matrix $R$ corresponding to a given partition. Then the eigenvalues
%of $B$ interlace the eigenvalues of $R$.
%\end{lem}

%The following lemma is also useful To compare the $A_\alpha$-spectral radius of a graph and its subgraph, we use the following lemma.

\begin{lem}(See \cite{go})\label{lem-2}
If $M_1$ and $M_2$ are two nonnegative $n\times n$ matrices such that  $M_1-M_2$ is nonnegative, then
$$\rho(M_1)\geq \rho(M_2),$$
where $\rho(M_i)$ is the spectral radius of $M_i$ for $i=1,2$.
\end{lem}

% The following result due to Tutte \cite{tt} shows that the converse is true. It is a special case of the Tutte–Berge Formula \cite{bo}:
% For any graph $G$,
% $$\alpha'(G)=\frac{1}{2} \textrm{min}\{V(G)-(\circ(G-S)-|S|):S \subseteq V \}$$

\begin{lem}(Tutte's Theorem, \cite{bo,tt})\label{lem-3}
A graph $G$ has a perfect matching if and only if $$o(G-S)\leq |S|~\mbox{for all}~S\subseteq V(G),$$ where $o(G-S)$ is the number of odd components of $G-S$.
\end{lem}

\section{Proof of Theorem \ref{thm-2}}

In this section, we shall give the proof of Theorem \ref{thm-2}.

\renewcommand\proofname{\bf Proof of Theorem \ref{thm-2}}
\begin{proof}
 Clearly, we have  $\alpha'(K_1\nabla(K_{n-3}\cup 2K_1))=n/2-1$, and so $K_1\nabla(K_{n-3}\cup 2K_1)$ has no perfect matchings. Assume that $G$ is a connected graph of order $n$ without perfect matchings. To prove Theorem \ref{thm-2}, it suffices to show that  $\rho_\alpha(G)\leq \rho_\alpha(K_1\nabla(K_{n-3}\cup 2K_1))$ with equality holding if and only if $G=K_1\nabla(K_{n-3}\cup 2K_1)$.

Since $G$ is connected and has no perfect matchings, by Lemma \ref{lem-3}, there exists some nonempty subset $S$ of $V (G)$ such that $o(G-S)-|S|\geq 1$ and all
components of $G-S$ are odd. Let $s=|S|$ and $q=o(G-S)$. Notice that $q$ and
$s$ must be of the same parity because $n$ is even.  Thus we have $q\geq s+ 2$.    Let $n_1, n_2, \ldots, n_q$ be the orders of the $q$ components of $G-S$, respectively. Without loss of generality, we take $n_1\geq n_2\geq \cdots\geq n_q$. Then $G$ is a spanning subgraph of $K_s\nabla(K_{n_1}\cup K_{n_2}\cup\cdots\cup K_{n_q})$, and so   $\rho_\alpha(G) \leq \rho_\alpha(K_s\nabla(K_{n_1}\cup K_{n_2}\cup\cdots\cup K_{n_q}))$ by Lemma \ref{lem-2}. We have the following three claims.
\begin{claim}\label{claim-1}
For $\alpha\in[0,1)$, we have
$$\rho_\alpha(K_s\nabla(K_{n_1}\cup K_{n_2}\cup\cdots\cup K_{n_q}))\leq \rho_\alpha(K_s\nabla(K_{n-s-q+1}\cup (q-1)K_1)),$$ where the equality holds if and only if $(n_1,n_2,\ldots,n_q)=(n-s-q+1,1,\ldots,1)$.
\end{claim}
\renewcommand\proofname{\bf Proof of Claim \ref{claim-1}}
\begin{proof}
Since $n_i$ is odd for all $i\in \{1,\ldots,q\}$,  it suffices to prove
$$\rho_\alpha(K_s\nabla(K_{n_1}\cup K_{n_2}\cup \cdots\cup K_{n_q})< \rho_\alpha(K_s\nabla(K_{n_1+2}\cup K_{n_2}\cup\cdots\cup K_{n_j-2}\cup \cdots\cup  K_{n_q})$$ whenever $n_j\geq 3$ for  $j\in\{2,\ldots,q\}$. Without loss of generality, we take $j=q$, and for other cases the proof is similar.

Let  $G_1=K_s\nabla(K_{n_1}\cup K_{n_2}\cup\cdots\cup K_{n_q})$ and $G_2=K_s\nabla(K_{n_1+2}\cup K_{n_2}\cup \cdots\cup K_{n_q-2})$. It is easy to see that $A_\alpha(G_1)$ has the  equitable quotient matrix
$$
B_1=
\begin{bmatrix}
n\alpha-s\alpha+s-1 & n_1(1-\alpha) & n_2(1-\alpha) & \cdots & n_q(1-\alpha) \\
s(1-\alpha) & s\alpha+n_1-1 & 0 & \cdots & 0 \\
s(1-\alpha) & 0 & s\alpha+n_2-1 & \cdots & 0 \\
\vdots & \vdots & \ddots & \vdots & \vdots \\
s(1-\alpha) & 0 & \cdots & 0 & s\alpha+n_q-1
\end{bmatrix}.
$$
By a simple calculation, we find that the characteristic polynomial of $B_1$ is equal to
\begin{equation}\label{equ-1}
\begin{aligned}
\varphi_{B_1}(x)&=(x-n\alpha+s\alpha-s+1)\cdot \prod_{j=1}^q(x-s\alpha-n_j+1)\\
&~~~ +s(1-\alpha)^2\cdot \sum_{i=1}^q(-1)^in_i\prod_{j\neq i}(x-s\alpha-n_j+1)\\
&=\psi(s,n_1,n_2,\ldots,n_q;x).\\
\end{aligned}
\end{equation}
By Lemma \ref{lem-1},   $\rho_\alpha(G_1)$ coincides with the largest root of  $\varphi_{B_1}(x)=0$.  Then, by Lemma \ref{lem-2}, we obtain
\begin{equation}\label{equ-2}
\rho_\alpha(G_1)\geq n\alpha-s\alpha+s-1.
\end{equation}
Considering that  $A_\alpha(G_1)$ is irreducible and $G_1$ contains $K_{n_1+s}$ as a proper subgraph, we have
\begin{equation}\label{equ-3}
\rho_\alpha(G_1)>\rho_\alpha(K_{n_1+s})=n_1+s-1
\end{equation}
by the Perron-Frobenius theorem and the Rayleigh quotient. Notice that $A_\alpha(G_2)$ has the equitable quotient matrix
$$
B_2=
\begin{bmatrix}
n\alpha-s\alpha+s-1 & (n_1+2)(1-\alpha) & n_2(1-\alpha) & \cdots & (n_q-2)(1-\alpha) \\
s(1-\alpha) & s\alpha+n_1+1 & 0 & \cdots & 0 \\
s(1-\alpha) & 0 & s\alpha+n_2-1 & \cdots & 0 \\
\vdots & \vdots & \ddots & \vdots & \vdots \\
s(1-\alpha) & 0 & \cdots & 0 & s\alpha+n_q-3
\end{bmatrix}.
$$
The characteristic polynomials of $B_2$ is
$$
\varphi_{B_2}(x)=\psi(s,n_1+2,n_2,\ldots,n_q-2;x),
$$
where $\psi(\cdot)$ is defined  in (\ref{equ-1}). As above, we see that $\rho_\alpha(G_2)$ is equal to the largest root of $\varphi_{B_2}(x)=0$. By a simple computation, we obtain
\begin{equation}\label{equ-4}
\begin{aligned}
\varphi_{B_2}(x)-\varphi_{B_1}(x)&=\psi(s,n_1+2,n_2,\ldots,n_q-2;x)-\psi(s,n_1,n_2,\ldots,n_q;x)\\
&=b_0(x)+2s(1-\alpha)^2\sum_{i=1}^qb_i(x),
\end{aligned}
\end{equation}
where
$$
\begin{aligned}
b_0(x)&=-(n_1-n_q+4)(x-n\alpha+s\alpha-s+1)\prod_{j=2}^{q-1}(x-s\alpha-n_j+1),\\
b_1(x)&=-(x+n_1-n_q-s\alpha +3)\prod_{j=2}^{q-1}(x-s\alpha-n_j+1),\\
b_i(x)&=(-1)^{i+1}n_i(n_1-n_q+4)\prod_{j=2,j\neq i}^{q-1}(x-s\alpha-n_j+1)~\mbox{for $i=2,\ldots,q-1$},\\
b_q(x)&=(-1)^{q+1}(x-n_1+n_q-s\alpha-1)\prod_{j=2}^{q-1}(x-s\alpha-n_j+1).
\end{aligned}
$$
Combining (\ref{equ-2}), (\ref{equ-3}) and $n_1\geq n_2\geq \cdots\geq n_q$, we get $b_0(\rho_\alpha(G_1))\leq 0$. If $q$ is odd,  again by (\ref{equ-3}), we obtian
$$
\begin{aligned}
&b_1(\rho_\alpha(G_1))+b_q(\rho_\alpha(G_1))=-2(n_1-n_q+4)\prod_{j=2}^{q-1}(\rho_\alpha(G_1)-s\alpha-n_j+1)<0,\\
&b_{q-1}(\rho_\alpha(G_1))=-n_{q-1}(n_1-n_q+4)\prod_{j=2}^{q-2}(\rho_\alpha(G_1)-s\alpha-n_j+1)<0,
\end{aligned}
$$
and
$$
\begin{aligned}
b_i(\rho_\alpha(G_1))+b_{i+1}(\rho_\alpha(G_1))&=-(n_i-n_{i+1})(n_1-n_q+4)(\rho_\alpha(G_1)-s\alpha +1)\cdot \\
&~~~~\prod_{j=2,j\neq i,i+1}^{q-1}(\rho_\alpha(G_1)-s\alpha-n_j+1)\leq 0,
\end{aligned}
$$
for all even $i$ with $2\leq i\leq q-3$. Similarly, if $q$ is even, we have $b_1(\rho_\alpha(G_1))+b_q(\rho_\alpha(G_1))<0$, and  $b_i(\rho_\alpha(G_1))+b_{i+1}(\rho_\alpha(G_1))\leq 0$ for all even $i$ with $2\leq i\leq q-2$. Since $\varphi_{B_1}(\rho_\alpha(G_1))=0$, in both cases,  from  the above arguments and (\ref{equ-4})  we deduce that
$$
\varphi_{B_2}(\rho_\alpha(G_1))=\varphi_{B_2}(\rho_\alpha(G_1))-\varphi_{B_1}(\rho_\alpha(G_1))<0,
$$
which implies that $\rho_\alpha(G_1)<\rho_\alpha(G_2)$. This proves Claim \ref{claim-1}.
\end{proof}

\begin{claim}\label{claim-2}
 For $\alpha\in [0,1)$, we have
 $$\rho_\alpha(K_s\nabla(K_{n-s-q+1}\cup (q-1)K_1))\leq \rho_\alpha(K_s\nabla(K_{n-2s-1}\cup (s+1)K_1)),$$
where the equality holds if and only if $q=s+2$.
 \end{claim}
\renewcommand\proofname{\bf Proof of Claim \ref{claim-2}}
\begin{proof}
Recall that $q\geq s+2$. It suffices to prove
$$\rho_\alpha(K_s\nabla(K_{n-s-q+1}\cup (q-1)K_1))< \rho_\alpha(K_s\nabla(K_{n-s-q+3}\cup (q-3)K_1))$$
for $q\geq s+4$. Let $G_3=K_s\nabla(K_{n-s-q+1}\cup (q-1)K_1)$ and $G_4=K_s\nabla(K_{n-s-q+3}\cup (q-3)K_1)$. Notice that $A_\alpha(G_3)$ and $A_\alpha(G_4)$ have the equitable quotient matrices
$$
B_3=\begin{bmatrix}
n\alpha-s\alpha+s-1  & (n-q-s+1)(1-\alpha) & (q-1)(1-\alpha) \\
s(1-\alpha) & n-q+s\alpha-s &0 \\
s(1-\alpha) &  0 & s\alpha
\end{bmatrix}
$$
and
$$
B_4=\begin{bmatrix}
n\alpha-s\alpha+s-1  & (n-q-s+3)(1-\alpha) & (q-3)(1-\alpha) \\
s(1-\alpha) & n-q+s\alpha-s+2 &0 \\
s(1-\alpha) &  0 & s\alpha
\end{bmatrix},
$$
for which we denote the characteristic polynomials  by $\varphi_{B_3}(x)$ and $\varphi_{B_4}(x)$, respectively. From $\varphi_{B_4}(\rho_\alpha(G_4))=0$ we obtain
$$
\begin{aligned}
&(\rho_\alpha(G_4)-n+q-s\alpha+s-2)[(\rho_\alpha(G_4)-s\alpha)(\rho_\alpha(G_4)-n\alpha+s\alpha-s+1)\\
&-s(1-\alpha)^2(q-3)]-(n-s-q+3)s(1-\alpha)^2(\rho_\alpha(G_4)-s\alpha)=0.
\end{aligned}
$$
Combining this with  $\rho_\alpha(G_4)>n-q+2>n-q+s\alpha-s+2>s\alpha$ (by (\ref{equ-3})), we have
\begin{equation}\label{equ-4-1}
(\rho_\alpha(G_4)-s\alpha)(\rho_\alpha(G_4)-n\alpha+s\alpha-s+1)-s(1-\alpha)^2(q-3)>0.
\end{equation}
On the other hand, a simple calculation yields that
$$
\begin{aligned}
\varphi_{B_3}(x)-\varphi_{B_4}(x)&=2s(1-\alpha)^2(n-q-s)\\
&~~~+2[(x-s\alpha)(x-n\alpha+s\alpha-s+1)-s(1-\alpha)^2(q-3)].
\end{aligned}
$$
Then from (\ref{equ-4-1}) we have
\begin{equation}\label{equ-5}
\begin{aligned}
\varphi_{B_3}(\rho_\alpha(G_4))&=\varphi_{B_3}(\rho_\alpha(G_4))-\varphi_{B_4}(\rho_\alpha(G_4))\\
&=2s(1-\alpha)^2(n-q-s)\\
&~~~+2\big[(\rho_\alpha(G_4)\!-\!s\alpha)(\rho_\alpha(G_4)\!-\!n\alpha\!+\!s\alpha\!-\!s\!+\!1)\!-\!s(1\!-\!\alpha)^2(q\!-\!3)\big]\\
&>0.
\end{aligned}
\end{equation}
Let $\gamma_1=\rho_\alpha(G_3)\geq \gamma_2\geq \gamma_3$ be the three roots of $\varphi_{B_3}(x)=0$. We claim that $\gamma_2<\rho_\alpha(G_4)$. In fact, let $D=\mathrm{diag}(s,n-q-s+1,q-1)$. It is easy to check that $D^{1/2}B_3D^{-1/2}$ is  symmetric, and also contains
$$
\begin{bmatrix} n-q+s\alpha-s & 0\\ 0 & s\alpha\end{bmatrix}
$$
as its submatrix. Since $D^{1/2}B_3D^{-1/2}$ and $B_3$ have the same eigenvalues,  the Cauchy interlacing theorem  (cf. \cite{hh}) implies that $\gamma_2\leq n-q+s\alpha-s< n-q+2<\rho_\alpha(G_4)$.  Therefore,  we conclude that $\rho_\alpha(G_3)<\rho_\alpha(G_4)$ by  (\ref{equ-5}), as required.
\end{proof}

%\begin{claim}\label{claim-3}
% If  $\alpha\in[0,2/3]$ and $n\geq 16$ or $\alpha\in(2/3,1)$ and $n\geq 5/(1-\alpha)$, we have
%$$\rho_\alpha(K_s\nabla(K_{n-2s-1}\cup (s+1)K_1))\leq \rho_\alpha(K_1\nabla(K_{n-3}\cup 2K_1)),$$
%with equality holding if and only if $s=1$.
%\end{claim}

\begin{claim}\label{claim-3}
 If  $\alpha\in[0,1)$ and $n\geq f(\alpha)$, where $f(\alpha)$ is defined in (\ref{equ-0}), then we have
$$\rho_\alpha(K_s\nabla(K_{n-2s-1}\cup (s+1)K_1))\leq \rho_\alpha(K_1\nabla(K_{n-3}\cup 2K_1)),$$
with equality holding if and only if $s=1$.
\end{claim}
\renewcommand\proofname{\bf Proof of Claim \ref{claim-3}}
\begin{proof}
Let $G_5^s=K_s\nabla(K_{n-2s-1}\cup (s+1)K_1)$. Clearly, $1\leq s\leq n/2-1$. We see that $A_\alpha(G_5^s)$ has the equitable quotient matrix
 $$
 B_5^s=
 \begin{bmatrix}
 n\alpha-s\alpha+s-1 & (n-2s-1)(1-\alpha) & (s+1)(1-\alpha)\\
 s(1-\alpha) & n+s\alpha-2s-2& 0\\
 s(1-\alpha) & 0 &s\alpha
 \end{bmatrix}.
 $$
By a simple computation, the characteristic polynomial of $B_5^s$ is equal to
\begin{equation}\label{equ-6}
\begin{aligned}
\varphi_{B_5^s}(x)&=x^3 + (3 \!-\! \alpha n \!-\! s (\alpha \!-\! 1) \!-\! n) x^2 - (s^2 \!-\! \alpha (\alpha n \!-\! 2) s \!-\! (\alpha n \!-\! 1) (n \!-\! 2)) x\\
 &~~~+(3\alpha-2)(1-\alpha)s^3 + (8 \alpha + n - 2 \alpha n + 2 \alpha^2 n - 5 \alpha^2 - 4) s^2\\
 &~~~ - (n - 2) (\alpha + \alpha^2 n - \alpha^2 - 1) s.
\end{aligned}
\end{equation}
Let $\gamma_1=\rho_\alpha(G_5^s)\geq \gamma_2\geq\gamma_3$ be the three roots of $\varphi_{B_5^s}(x)=0$. As in Claim \ref{claim-2}, by using the Cauchy interlacing theorem, we can deduce that $\gamma_2\leq n+s\alpha-2s-2<n-3 $. Also, since $\varphi_{B_5^1}(n-3)=2(\alpha - 1)((n-5)\alpha+ 1)<0$, we have  $\rho_\alpha(G_5^1)>n-3$.  We consider the following two situations.

{\flushleft \bf Case 1.} $0\leq \alpha\leq 2/3$ and $n\geq f(\alpha)$.

First we assert that $\rho_\alpha(G_5^{n/2-1})<\rho_\alpha(G_5^{1})$. In fact, by calculating the largest eigenvalue of $B_5^{n/2-1}$, we obtain
$$
\rho_\alpha(G_5^{n/2-1})=\frac{1}{4}\big((2\alpha+1)n-4+\sqrt{(4\alpha^2 - 8\alpha + 5)n^2 + 8(\alpha - 1)n}\big).
$$
For $0\leq \alpha\leq 1/2$ and $n\geq 12$, or $1/2<\alpha\leq 2/3$ and $n\geq 16$, we have 
$$
(1-\alpha)n^2 + (6\alpha-10)n + 16 >0,
$$
which implies that $\rho_\alpha(G_5^{n/2-1})< n-3<\rho_\alpha(G_5^{1})$, as required. For   $0\leq \alpha\leq 1/2$ and $n=10$, by a simple computation, we have
$$
\varphi_{B_5^1}(\rho_\alpha(G_5^4))=-\frac{3}{2}\left[50\alpha^2-99\alpha+34+(5\alpha-2)\sqrt{100\alpha^2-180\alpha+105}\right]<0,
$$
and so $\rho_\alpha(G_5^{n/2-1})=\rho_\alpha(G_5^{4})<\rho_\alpha(G_5^{1})$. Similarly, for $1/2<\alpha\leq 2/3$ and $n=14$, we obtain 
$$\varphi_{B_5^1}(\rho_\alpha(G_5^6))=-\frac{5}{2}\left[102\alpha^2 - 215\alpha + 78+(7\alpha-2)\sqrt{196\alpha^2 - 364\alpha + 217}\right]<0,$$ implying that $\rho_\alpha(G_5^{n/2-1})=\rho_\alpha(G_5^{6})<\rho_\alpha(G_5^{1})$. Therefore, we conclude that  $\rho_\alpha(G_5^{n/2-1})<\rho_\alpha(G_5^{1})$ for $0\leq \alpha\leq 2/3$ and $n\geq f(\alpha)$.

Now suppose $2\leq s\leq n/2-2$. A simple calculation yields that
$$
\begin{aligned}
\varphi_{B_5^s}(x)-\varphi_{B_5^1}(x)&=(s-1)[(1-\alpha)x^2+ (n\alpha^2 - 2\alpha - s - 1)x\\
&~~~~ -\alpha^2 n^2+ ((2s + 5)\alpha^2 - ( 2s + 3)\alpha + s + 2)n \\
&~~~~- (3s^2 + 8s + 10)\alpha^2 + ( 5s^2 + 13s + 15)\alpha - 2s^2 - 6s - 8].\\
\end{aligned}
$$
Since $n\geq 2s+4$ and $s\geq 2$, we have
$$-(n\alpha^2 - 2\alpha - s - 1)/(2(1-\alpha))<n-3<\rho_\alpha(G_5^1).$$ It follows that
\begin{equation}\label{equ-7}
\begin{aligned}
\varphi_{B_5^s}(\rho_\alpha(G_5^1))&=\varphi_{B_5^s}(\rho_\alpha(G_5^1))-\varphi_{B_5^1}(\rho_\alpha(G_5^1))\\
&>(s-1)[(1-\alpha)(n-3)^2+ (n\alpha^2 - 2\alpha - s - 1)(n-3)\\
&~~~ -\alpha^2 n^2+ ((2s + 5)\alpha^2 - ( 2s + 3)\alpha + s + 2)n \\
&~~~- (3s^2 + 8s + 10)\alpha^2 + ( 5s^2 + 13s + 15)\alpha - 2s^2 - 6s - 8]\\
&=(s-1)[(1-\alpha)n^2+ ((2s + 2)\alpha^2 - (2s-1)\alpha - 5)n\\
&~~~~-(3s^2 + 8s + 10)\alpha^2 + (5s^2 + 13s + 12)\alpha - 2s^2 - 3s + 4].\\
\end{aligned}
\end{equation}
If $s\geq 4$, one can verify that $-((2s + 2)\alpha^2 - (2s-1)\alpha - 5)/(2(1-\alpha))<2s+4$. From (\ref{equ-7}) we have
$$
\begin{aligned}
\varphi_{B_5^s}(\rho_\alpha(G_5^1)) &>(s-1)[(1-\alpha)(2s+4)^2+ ((2s + 2)\alpha^2 - (2s-1)\alpha - 5)(2s+4)\\
&~~~~-(3s^2 + 8s + 10)\alpha^2 + (5s^2 + 13s + 12)\alpha - 2s^2 - 3s + 4]\\
&=(s-1)[(s^2 + 4s - 2)\alpha^2 -(3s^2 + 9s)\alpha + 2s^2 + 3s]\\
&\geq (s-1)\left[\frac{4}{9}(s^2 + 4s - 2)  -\frac{2}{3}(3s^2 + 9s) + 2s^2 + 3s\right]\\
&=\frac{1}{9}(s-1)(4s^2 - 11s - 8)\\
&> 0,
\end{aligned}
$$
where the last two inequalities follow from $(3s^2 + 9s)/(2(s^2 + 4s - 2))>2/3\geq \alpha$ and $s\geq 4$, respectively. If $s=2$, then  $-((2s + 2)\alpha^2 - (2s-1)\alpha - 5)/(2(1-\alpha))=-(6\alpha^2-3\alpha-5)/(2(1-\alpha))<10\leq f(\alpha)\leq n$ due to $0\leq \alpha\leq 2/3$. Again by (\ref{equ-7}),  we obtain
$$
\begin{aligned}
\varphi_{B_5^s}(\rho_\alpha(G_5^1)) &>(s-1)[10^2\cdot (1-\alpha)+ 10\cdot ((2s + 2)\alpha^2 - (2s-1)\alpha - 5)\\
&~~~~-(3s^2 + 8s + 10)\alpha^2 + (5s^2 + 13s + 12)\alpha - 2s^2 - 3s + 4]\\
&=22\alpha^2 - 72\alpha + 40\\
&> 0.
\end{aligned}
$$
Similarly, if $s=3$, then $-((2s + 2)\alpha^2 - (2s-1)\alpha - 5)/(2(1-\alpha))=-(8\alpha^2-5\alpha-5)/(2(1-\alpha))<10\leq f(\alpha)\leq n$ due to $0\leq \alpha\leq 2/3$. Thus, for $0\leq \alpha\leq 1/2$ and $n\geq f(\alpha)=10$, we have
$$
\begin{aligned}
\varphi_{B_5^s}(\rho_\alpha(G_5^1)) &>(s-1)[10^2\cdot (1-\alpha)+ 10\cdot ((2s + 2)\alpha^2 - (2s-1)\alpha - 5)\\
&~~~~-(3s^2 + 8s + 10)\alpha^2 + (5s^2 + 13s + 12)\alpha - 2s^2 - 3s + 4]\\
&=38\alpha^2 - 108\alpha + 54\\
&> 0,
\end{aligned}
$$
and for $1/2<\alpha\leq 2/3$ and $n\geq f(\alpha)=14$,  we have
$$
\begin{aligned}
\varphi_{B_5^s}(\rho_\alpha(G_5^1)) &>(s-1)[14^2\cdot (1-\alpha)+ 14\cdot ((2s + 2)\alpha^2 - (2s-1)\alpha - 5)\\
&~~~~-(3s^2 + 8s + 10)\alpha^2 + (5s^2 + 13s + 12)\alpha - 2s^2 - 3s + 4]\\
&=102\alpha^2 - 340\alpha + 206\\
&> 0.
\end{aligned}
$$
Therefore, we conclude that $\varphi_{B_5^s}(\rho_\alpha(G_5^1))>0$ for $2\leq s\leq n/2-2$.
As  $\gamma_2<n-3<\rho_\alpha(G_5^1)$, we  have $\rho_\alpha(G_5^s)<\rho_\alpha(G_5^1)$ for $2\leq s\leq n/2-2$, and so the result follows.

{\flushleft \bf Case 2.} $2/3<\alpha<1$ and $n\geq f(\alpha)=5/(1-\alpha)$.

In this situation, we shall prove that  $\rho_\alpha(G_5^s)\leq n-3<\rho_\alpha(G_5^1)$ for $2\leq s\leq n/2-1$. As $\gamma_2<n-3$, it suffices to show that $\varphi_{B_5^s}(n-3)>0$ for $2\leq s\leq n/2-1$. According to (\ref{equ-6}), we have
$$
\begin{aligned}
\varphi_{B_5^s}(n-3)&=(3\alpha-2)(1-\alpha)s^3 + ((2 \alpha^2 - 2 \alpha) n - 5 \alpha^2 + 8 \alpha - 1) s^2\\
&~~~ + ((1 - \alpha) n^2 + (3 \alpha - 5) n - 2 \alpha^2 - \alpha + 7) s + (\alpha - 1)n^2 + (5 - 3\alpha)n - 6\\
&:=\psi(s,n).
\end{aligned}
$$
Then we see that
$$
\begin{aligned}
\frac{\partial\psi(s,n)}{\partial s} &=3(3\alpha-2)(1-\alpha)s^2 + 2((2 \alpha^2 - 2 \alpha) n - 5 \alpha^2 + 8 \alpha - 1) s\\
&~~~+ ((1 - \alpha) n^2 + (3 \alpha - 5) n - 2 \alpha^2 - \alpha + 7).
\end{aligned}
$$
Since $2/3<\alpha<1$ and $n\geq 5/(1-\alpha)$, by a simple calculation, we have
$$
\begin{aligned}
\frac{\partial\psi(s,n)}{\partial s}\mid_{s=2}&=(1 - \alpha) n^2 + (8 \alpha^2 - 5 \alpha - 5) n - 58 \alpha^2 + 91 \alpha - 21\\
&\geq (1 - \alpha) \Big(\frac{5}{1-\alpha}\Big)^2 + (8 \alpha^2 - 5 \alpha - 5) \Big(\frac{5}{1-\alpha}\Big) - 58 \alpha^2 + 91 \alpha - 21\\
&=\frac{1}{1-\alpha}(58\alpha^3 - 109\alpha^2 + 87\alpha - 21)\\
&>0,
\end{aligned}
$$
and
$$
\begin{aligned}
\frac{\partial\psi(s,n)}{\partial s}\mid_{s=\frac{n}{2}-1}&=\frac{1}{4}(1-\alpha)((\alpha - 2)n^2 + 4\alpha + 12)\\
&\leq \frac{1}{4}(1-\alpha)\Big((\alpha - 2)\Big(\frac{5}{1-\alpha}\Big)^2 + 4\alpha + 12\Big)\\
&= \frac{1}{4(1-\alpha)}(4\alpha^3 + 4\alpha^2 + 5\alpha- 38)\\
&<0.
\end{aligned}
$$
This implies that $\varphi_{B_5^s}(n-3)=\psi(s,n)\geq \min\{\psi(2,n),\psi(n/2-1,n)\}$ because  the leading coefficient of $\psi(s,n)$ (viewed as a cubic polynomial of $s$) is positive, and $2\leq s\leq n/2-1$. Again by $2/3<\alpha<1$ and $n\geq 5/(1-\alpha)$, we have
$$
\begin{aligned}
\psi(2,n)&=(1 - \alpha)n^2 + (8\alpha^2 - 5\alpha - 5)n - 48\alpha^2 + 70\alpha - 12\\
&\geq(1 - \alpha)\Big(\frac{5}{1-\alpha}\Big)^2 + (8\alpha^2 - 5\alpha - 5)\Big(\frac{5}{1-\alpha}\Big) - 48\alpha^2 + 70\alpha - 12\\
&=\frac{3}{(1-\alpha)}(16\alpha^3 - 26\alpha^2 + 19\alpha - 4)\\
&>0,
\end{aligned}
$$
and
$$
\begin{aligned}
\psi(n/2-1,n)&=\frac{1}{8}((2-\alpha)n + 2\alpha-6)((1-\alpha)n^2 + (6\alpha-10)n + 16)\\
&\geq \frac{1}{8}\Big((2\!-\!\alpha)\Big(\frac{5}{1\!-\!\alpha}\Big) \!+\! 2\alpha\!-\!6\Big)\Big((1\!-\!\alpha)\Big(\frac{5}{1\!-\!\alpha}\Big)^2 \!+\! (6\alpha\!-\!10)\Big(\frac{5}{1\!-\!\alpha}\Big) \!+\! 16\Big)\\
&=\frac{1}{8(1-\alpha)^2}(- 28\alpha^3 + 60\alpha^2 + 29\alpha - 36)\\
&>0.
\end{aligned}
$$
Therefore, for $2\leq s\leq n/2-1$, we conclude that $\varphi_{B_5^s}(n-3)\geq \min\{\psi(2,n),\psi(n/2-1,n)\}>0$, as required. This proves Claim \ref{claim-3}.
\end{proof}

From Claims \ref{claim-1}--\ref{claim-3} and the arguments at the beginning  of the proof, we conclude that $\rho_\alpha(G)\leq \rho_\alpha(K_1\nabla(K_{n-3}\cup 2K_1))$, with equality holding if and only if $G=K_1\nabla(K_{n-3}\cup 2K_1)$. Also note that $\rho_\alpha(K_1\nabla(K_{n-3}\cup 2K_1))$ is the largest root of $\varphi_{B_5^1}(x)=x^3 - ((\alpha + 1)n +\alpha-4)x^2 + (\alpha n^2 + (\alpha^2 - 2\alpha - 1)n - 2\alpha+1)x -\alpha^2n^2 + (5\alpha^2 - 3\alpha + 2)n - 10\alpha^2 + 15\alpha - 8=0$.

We complete the proof.
\end{proof}

\begin{remark}
\emph{
In Theorem \ref{thm-2}, we determine the unique extremal graph for  each $n\geq f(\alpha)$, where $\alpha\in[0,1)$ and $f(\alpha)$ is defined in (\ref{equ-0}). For the remaining (not too many) cases, i.e., $n<f(\alpha)$, according to the proof of Theorem \ref{thm-2}, we know that if $\rho_\alpha(G)>\max\{\rho_\alpha(K_s\nabla(K_{n-2s-1}\cup (s+1)K_1)):1\leq s\leq n/2-1\}$ then $G$ has a perfect matching. Thus it remains to determine the maximum value of $\rho_\alpha(K_s\nabla(K_{n-2s-1}\cup (s+1)K_1))$ when $s$ ranges over  $[1,n/2-1]$. However, there are two aspects of difficulties. The first one is that $\rho_\alpha(K_s\nabla(K_{n-2s-1}\cup (s+1)K_1))$ is the root of a (usually not factorable) cubic polynomial whose coefficients are dependent on $\alpha$ and $n$ (see (\ref{equ-6})), and it is usually hard to compare these roots. The second one is that the maximum value might not be attained at $s=1$. For example, for $\alpha=2/3$ and $n=12$, one can verify that the maximum value is attained at $s=n/2-1=5$, which suggests that the bound $n\geq 14$ for $\alpha=2/3$ in Theorem \ref{thm-2} cannot be improved. Also, we can check that,  for $(\alpha,n)=(1/4,8)$, $(11/16,12)$, $(3/4,16)$ or $(9/10,40)$, the maximum value is attained at $s=n/2-1$, but for  $(\alpha,n)=(1/8,8)$, $(11/16,14)$, $(3/4,18)$ or $(9/10,42)$, the maximum value is attained at $s=1$. Indeed, by examining many special values of $\alpha$ and $n$, we find that the maximum value is always attained at $s=1$ or $s=n/2-1$.
}
\end{remark}

\section*{Acknowledgements}
The authors are indebted to the anonymous referee for his/her  valuable comments and helpful suggestions. X. Huang is partially supported  by the National Natural Science Foundation of China (Grant No. 11901540 and Grant No. 11671344).

\end{document}